\documentclass[11pt,a4paper]{article}
\usepackage[english]{babel}
\usepackage{amsmath, amsthm, amsfonts, amssymb, amscd}
\usepackage[utf8]{inputenc}
\usepackage[pdftex]{graphicx}
\usepackage[matrix,arrow]{xy}
\usepackage{enumerate}
\DeclareGraphicsExtensions{.bmp,.png,.pdf,.jpg}
\newtheorem*{proper}{Property}
\newtheorem{defi}{Definition}

\newtheorem{teo}{Theorem}
\newtheorem{lema}[teo]{Lemma}
\newtheorem{cor}[teo]{Corollary}
\newtheorem{prop}[teo]{Proposition}

\newcommand{\hookuparrow}{\mathrel{\rotatebox[origin=c]{90}{$\hookrightarrow$}}}

\newcommand{\pf}{\noindent{\bf Proof\ \ }}
\newcommand{\cqd}{{\hfill $\rule{2mm}{2mm}$}\vspace{3mm}}
\newcommand{\Z}{\mathbb Z}
\newcommand{\N}{\mathbb N}

\newcommand{\I}{\mathcal{I}}
\newcommand{\J}{\mathcal{J}}

\newcommand{\R}{\mathcal{R}}

\usepackage{float}

\begin{document}
\title{On the Colength of Fractional Ideals}

\author{E. M. N. de Guzm\'an  
\and A. Hefez} 

\date{\small{Departamento de Matem\'atica Aplicada \\ Universidade Federal Fluminense }}
\maketitle

\noindent{\bf Abstract}
The main result in this paper is to supply a recursive formula, on the number of minimal primes, for the colength of a fractional ideal in terms of the maximal points of the value set of the ideal itself. The fractional ideals are taken in the class of complete admissible rings, a more general class of rings than those of algebroid curves. For such rings with two or three minimal primes, a closed formula for that colength is provided, so improving results by Barucci, D'Anna and Fr\"oberg.\bigskip

\medskip

\section{Introduction}

The computation of the colength of a fractional ideal of a ring of an algebroid plane branch in terms of its value set was known at least since the work of Gorenstein in the fifties of last century (cf. \cite{Go}). Another instance of such computation was performed for a larger class of analytically reduced but reducible rings by D'Anna (cf. \cite[\S 2]{Da}), who related such a colength to the length of a maximal saturated chain in the set of values of the given fractional ideal. The issue with D'Anna's method is that it requires the knowledge of many elements in the set of values, a feature that would be desirable to overcome to increase the computational effectivity. In fact, in the particular case of an algebroid curve with two branches, Barucci, D'Anna and Fr\"oberg, in \cite{BDF1}, were able to give an explicit formula for the colength of a given fractional ideal in terms of fewer points of its value set, namely, those enjoying a certain maximal property (maximal points). Local rings of algebroid curves, as well as the larger class studied by D'Anna in \cite{Da} are a particular case of a larger class of local rings we shall be concerned with in this paper and that we refer to as {\em admissible}. By an admissible ring we shall mean a one dimensional, local, noetherian, Cohen-Macaulay, analytically reduced and residually rational ring such that the cardinality of its residue field is sufficiently large (see \cite{KST} for more details). For simplicity and without loss of generality (cf. \cite[\S 1]{Da}), we will also assume that our rings are complete with respect to the topology induced by the maximal ideal. In such case, a sufficiently large residue field means that its cardinality is greater than or equal to the number of its minimal primes $r$. One of the main results of this paper is Theorem 10. It gives a recursive formula on the number $r$, that corresponds to the number of branches in the case of rings of algebroid curves, for the colength of a fractional ideal in a complete admissible ring. The important feature consists that the computation requires only a few points of the value set (some special maximal ones). The other main result is Corollary 19 that provides a closed formula for the colength in the case of three minimal primes. It is worth to notice that such a closed formula for three minimal primes is not exactly a straightforward consequence of the recursive formula established in Theorem 10, as its proof demands a careful analysis of the geometry of the maximal points of the value set, with respect to the natural order relation in $\Z^r$ as recalled in Section 3. The outline of the paper is as follows. Section 2 collects some preliminaries and notation regarding the general background of the article. Section 3 is concerned with the definition of value sets and the partial order inherited by that of $\Z^r$, recalling three useful properties analogous to ones obtained for semigroups of values by Delgado and Garcia (cf. \cite{D87} and \cite{Ga}). Section 4 introduces and analyzes different kinds of maximal points in the value set to get enough tools to pass to Section 5 that is eventually concerned with the announced recursive formula for the colength of fractional ideals in admissible rings. To ease the comparison with the previous results due to Barucci, D'Anna and Fr\"oberg, we first analyze their recipe for $r=2$, while we devote Section 5.2 to the case $r\geq 3$: Basing on the Key Lemma 9, one eventually obtains, patching together all the pieces of the puzzle, the announced Theorem 10 that substantially  improves the nice method by D'Anna that, nevertheless, has the drawback that it needs the knowledge of many points in the value set. The closed formula for $r=3$ is finally dealt with in Section 6 where a fine detailed analysis of the geometry of the maximal points is offered in a series of lemmas, culminating with Lemma 17 that unavoidably leads, after the case by case analysis, the statement and proof of Theorem 18. It confirms a conjectural formula by M. Hernandes and implies Corollary 19, the closed formula for the colength in the case $r=3$ minimal primes that, in any case, furnish a substantial improvement of the results by Barucci, et al.\bigskip

\noindent {\bf Acknowledgements}	The first author was supported by a fellowship from CAPES, while the second one was partially supported by the CNPq Grant number 307873/2016-1. 
	
\section{General background}

In this section we refer to \cite{Da} for our unproved statements. Let $\wp_1, \ldots,\wp_r$ be the minimal primes of $R$. We will use the notation $I=\{1,\ldots,r\}$. We set $R_i=R/\wp_i$ and will denote by $\pi_i\colon R \to R_i$ the canonical surjection. 
Since $R$ is reduced, we have that $\bigcap_{i=1}^r \wp_i =\sqrt{(0)}=(0)$, so we get an injective homomorphism  
$$
\begin{array}{rcl}
\pi\colon R & \hookrightarrow &R_1\times \cdots \times R_r\\
            h & \mapsto & (\pi_1(h),\ldots,\pi_r(h)) .
						\end{array}$$

More generally, if $J=\{i_1,\ldots,i_s\}$ is any subset of $I$, we may consider $R_J=R/\cap_{j=1}^s \wp_{i_j}$ and will denote by $
\pi_J\colon R \longrightarrow R_J$ the natural surjection. 
\medskip

We will denote by $K$ the total ring of fractions of $R$ and when $J\subset I$ we denote by $K_J$ the total ring of fractions of the ring
$R_J$. Notice that $R_I=R$ and $K_I=K$. If $J=\{i\}$, then $R_{\{i\}}$ is equal to the above defined domain $R_i$ whose field of fractions will be denoted by $K_i$. Let $\widetilde{R}$ be the integral closure of $R$ in $K$ and $\widetilde{R}_J$ be that of $R_J$ in $K_J$. One has that $\widetilde{R}_J\simeq \widetilde{R}_{i_1} \times \cdots \times \widetilde{R}_{i_s}$, which in turn is the integral closure of $ {R}_{i_1} \times \cdots \times {R}_{i_s}$ in its total ring of fractions. \medskip

We have the following diagram:
\[
\begin{array}{rcc}
K_J & \simeq & K_{i_1} \times \cdots \times K_{i_s} \\
\hookuparrow     &  & \hookuparrow \\
\widetilde{R}_J &  \simeq & \widetilde{R}_{i_1} \times \cdots \times \widetilde{R}_{i_s}\\
\hookuparrow     &   & \hookuparrow\\
R_J         & \hookrightarrow & R_{i_1} \times \cdots \times R_{i_s}	
	\end{array}
	\]
		
Since each $\widetilde{R_i}$ is a DVR, with a valuation denoted by $v_i$, one has that $K_i$  is a valuated field with the extension of the valuation $v_i$ which is denoted by the same symbol. This allows one to define the value map
\[
\begin{array}{ccc}
v \colon K \setminus Z(K) & \to &  \Z^r\\
   h & \mapsto &  (v_1(\pi_1(h)), \ldots ,v_r(\pi_r(h))),
	\end{array}
\]
where $\pi_i$ here denotes the projection $K \to \mathcal K_i$, which is the extension of the previously defined projection map $\pi_i\colon R \to R_i$ and $Z(K)$ stands for the set of zero divisors of $K$. \medskip

An $R$-submodule $\mathcal I$ of $K$ will be called a {\em regular fractional ideal} of $R$ if it  contains a regular element of $R$ and there is a regular element $d$ in $R$ such that $d\,\mathcal I \subset R$.

Since $d\,\mathcal I$ is an ideal of $R$, which is a noetherian ring, one has that $\mathcal I\subset K$ is a nontrivial fractional ideal if and only if it contains a regular element of $R$ and it is a finitely generated $R$-module.

Examples of fractional ideals of $R$ are $R$ itself, $\widetilde{R}$, the conductor $\mathcal C$ of $\tilde{R}$ in $R$, or any ideal of $R$ or of $\tilde{R}$ that contains a regular element. Also, if $\I$ is a regular fractional ideal of $R$, then for all $\emptyset \neq J\subset I$ one has that $\pi_J(\I)$ is a regular fractional ideal of $R_J$, where, this time, $\pi_J\colon K\to K_J$ denotes the natural projection.\medskip

\section{Value sets}

If $\mathcal I$ is a regular fractional ideal of $R$, we define the {\em value set} of $\mathcal I$ as being 
$$E=v(\mathcal I\setminus Z(K))\subset \Z^r.$$
 If $J=\{i_1,\ldots,i_s\}\subset I$, then we denote by
 $pr_J$ the projection $\Z^r \to \Z^s$, $(\alpha_1,\ldots,\alpha_r) \mapsto (\alpha_{i_1},\ldots,\alpha_{i_s})$. \medskip

Let us define 
$$E_J=v(\pi_J(\I)\setminus Z(K_J)). 
$$

If $j\in J=\{i_1,\ldots,i_t,\ldots i_s\} \subset I$, with $i_t=j$, for $\beta=(\beta_{i_1},\ldots,\beta_{i_s})\in E_J$, then we define $pr_j(\beta)=\beta_{i_j}$ and $$\widetilde{pr}_j(\beta)=\beta_{i_t}=\beta_j.$$ 
Notice that if $J=I$, then $\widetilde{pr}_i=pr_i$, for all $i\in I$.\medskip

We will consider on $\Z^r$ the product order $\leq$ and will write $(a_1,\ldots,a_r) < (b_1,\ldots,b_r)$ when $a_i<b_i$, for all $ i=1,\ldots,r$.\medskip

Value sets of fractional ideals have the following fundamental properties, analogous to the properties of semigroups of values described by Garcia for $r=2$ in \cite{Ga} and by Delgado for $r> 2$ in \cite{D87} (see also \cite{Da} or \cite{BDF1}):

\begin{proper}[\textbf{A}]
	If $\alpha=(\alpha_1,\ldots,\alpha_r)$ and $\beta=(\beta_1,\ldots,\beta_r)$ belong to $E$, then 
	$$min(\alpha,\beta)=(min(\alpha_1,\beta_1),\ldots,min(\alpha_r,\beta_r))\in E.$$
\end{proper}

\begin{proper}[\textbf{B}]
	If $\alpha=(\alpha_1,\ldots,\alpha_r), \beta=(\beta_1,\ldots,\beta_r)$ belong to $E$, $\alpha\neq\beta$ and $\alpha_i=\beta_i$ for some $i\in\{1,\ldots,r\}$, then there exists $\gamma\in E$ such that $\gamma_i>\alpha_i=\beta_i$ and $\gamma_j\geq min\{\alpha_j,\beta_j\}$ for each $j\neq i$, with equality holding if $\alpha_j\neq\beta_j$.
\end{proper}

\begin{proper}[\textbf{C}]
There exist $\alpha \in \Z^r$ and $\gamma\in \N^r$ such that
\[
\gamma + \N^r \subset E \subset \alpha +\Z^r.
\]
\end{proper}

Properties (A) and (C) allow one to conclude that there exist a unique $m_E=(m_1,\ldots,m_r)$ such that $\beta_i\geq m_i$, $i=1,\ldots,r$, for all $(\beta_1,\ldots,\beta_r)\in E$ and a unique least element $\gamma\in E$ with the property that $\gamma+\N^r \subset E$. This element is what we call the conductor of $E$ and will denote it by $c(E)$.

Observe that one always has
\[
c(E_J)\leq pr_J(c(E)), \ \ \forall\, J\subset I.
\]

One has the following result:

\begin{lema} \label{proj} If $\I$ is a fractional ideal of $R$ and $\emptyset \neq J\subset I$, then
$pr_J(E)=E_J$.
\end{lema}
\pf One has obviously that $pr_J(E)\subset E_J$. On the other hand, let $\alpha_J \in E_J$. Take $h\in \I$ such that $v_J(\pi_J(h))=\alpha_J$. If $h\not\in Z(K)$ we are done. Otherwise, choose any $h'\in \I \setminus Z(K)$ such that $pr_J(v(h'))>\alpha_J$, which exists since $E$ has a conductor. Hence, $v_J(h+h')=\alpha_J$, proving the other inclusion.
\cqd

\section{Maximal points} \label{maximals}

We now introduce the important notion of a {\em fiber} of an element $\alpha\in E$ with respect to a subset $J\subset I=\{1,\ldots,r\}$ that will play a central role in what follows.

\begin{defi} 	Given $A\subset \mathbb{Z}^r$, $\alpha\in\mathbb{Z}^r$ and $\emptyset \neq J\subset I$, we define 	
$$F_J(A,\alpha)=\{\beta\in A; pr_J(\beta)=pr_J(\alpha) \ \text{and} \ pr_{I\setminus J}(\beta)>pr_{I\setminus J}(\alpha)\},$$
$$\overline{F}_J(A,\alpha)=\{\beta\in A; pr_J(\beta)=pr_J(\alpha), \  \text{and} \ pr_{I\setminus J}(\beta) \geq pr_{I\setminus J}(\alpha)\},$$
\begin{equation} \label{fiber} F(A,\alpha)=\bigcup_{i=1}^rF_{\{i\}}(A,\alpha).
\end{equation}
	
The set (\ref{fiber}) will be called the {\em fiber} of $\alpha$ in $A$.
	
\end{defi}

The sets $F_{\{i\}}(A,\alpha)$ and $\overline{F}_{\{i\}}(A,\alpha)$ will be denoted simply by $F_{i}(A,\alpha)$ and $\overline{F}_{i}(A,\alpha)$. Notice that $F_I(\mathbb{Z}^r,\alpha)=\overline{F}_I(\mathbb{Z}^r,\alpha)=\{\alpha\}$.

\begin{defi}
Let $\alpha\in A$. We will say that $\alpha$ is a\textbf{ maximal} point of $A$ if $F(A,\alpha)=\emptyset$. 
\end{defi}

This means that there is no element in $A$ with one coordinate equal to the corresponding coordinate of $\alpha$ and the others bigger.	\medskip

From now on $E$ will denote the value set of the regular fractional ideal $\I$ of $R$.

From the fact that $E$ has a minimum $m$ and a conductor $\gamma=c(E)$, one has immediately that all maximal elements of $E$ are in the limited region
\[
\{(x_1,\ldots,x_r)\in \Z^r; \ m_{i}\leq x_i < \gamma_i, \ \ i=1,\ldots,r\}.
\]

This implies that $E$ has finitely many maximal points.

\begin{defi} We will say that a maximal point $\alpha$ of $E$ is an \textbf{absolute maximal} if $F_J(E,\alpha)=\emptyset$ for every $J\subset I$, $J\neq I$. If a maximal point $\alpha$ of $E$ is such that $F_J(E,\alpha)\neq\emptyset$, for every $J\subset I$ with $\#J\geq2$, then $\alpha$ will be called a \textbf{relative maximal} of $E$.
\end{defi}

\begin{figure}[h]
 \centering
\includegraphics[scale=.6]{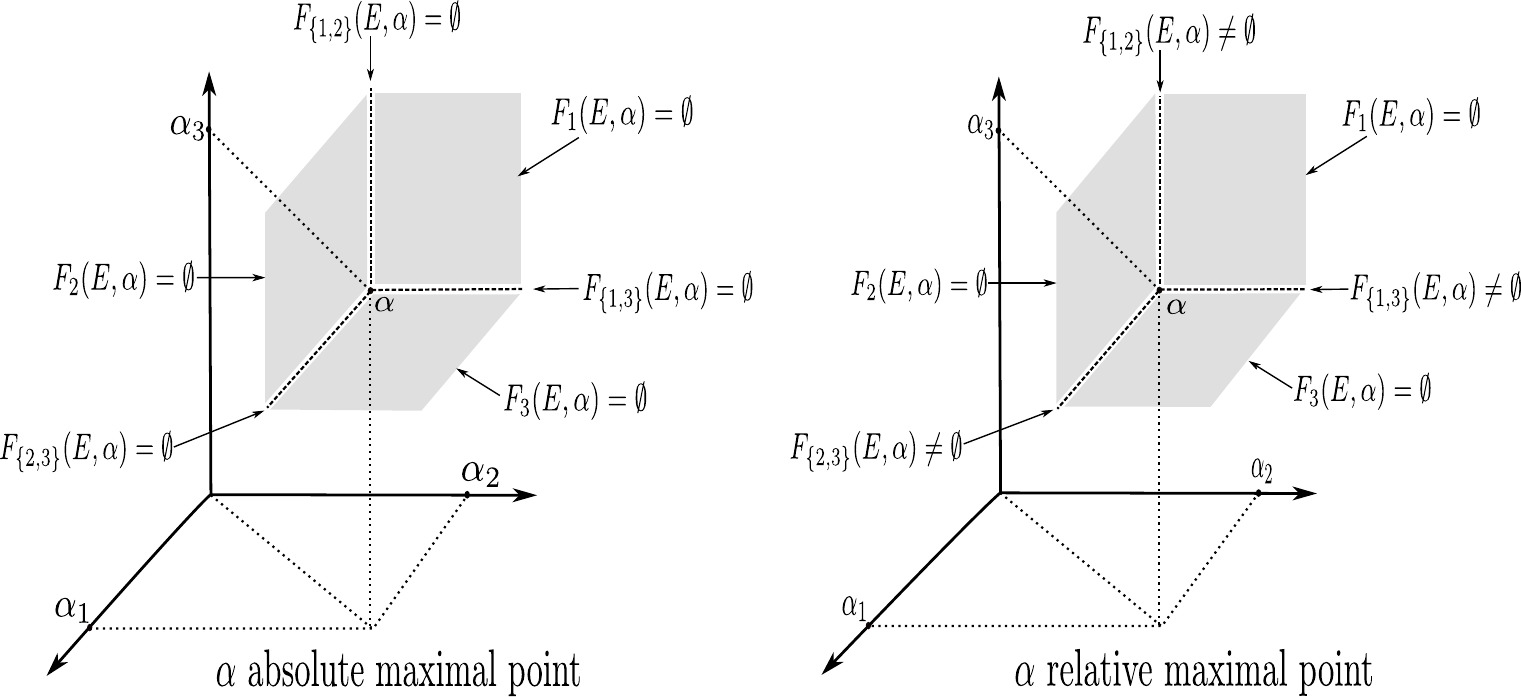}
\caption{Maximal points}
\end{figure}

In the case where $r=2$, the notions of maximal, relative maximal and absolute maximal coincide. For $r=3$ we may only have relative maximals or absolute maximals, but in general there will be several types of maximals.\medskip

We will denote by $M(E)$, $RM(E)$ and $AM(E)$ the sets of maximals, of relative maximals and absolute maximals of the set $E$, respectively.\medskip

The importance of the relative maximals is attested by the theorem below that says that the set $RM(E)$ determines $E$ in a combinatorial sense as follows:
\begin{teo}[generation]\label{generation}
	Let $\alpha\in\mathbb{Z}^r$ be such that $p_J(\alpha)\in E_J$ for all $J\subset I$ with $\#J=r-1$. Then
	$$\alpha\in E \Longleftrightarrow \alpha\notin F(\mathbb{Z}^r,\beta), \ \forall \beta\in RM(E).$$
\end{teo}

We will omit the proof since this result is a slight modification of \cite[Theorem 1.5 ]{D87} with essentially the same proof.\medskip

The following two lemmas give us characterizations of the relative and absolute maximal points that will be useful in Section 4.

\begin{lema}\label{carrel}
Given a value set $E\subset\mathbb{Z}^d$ and $\alpha\in \mathbb{Z}^d$ with the following properties:
\begin{enumerate}[\rm i)]
\item there is $i\in I$ such that $F_i(E,\alpha)=\emptyset$,
\item $F_{i,j}(E,\alpha)\neq\emptyset$ for all $j\in I\setminus\{i\}$.
\end{enumerate}
Then $\alpha$ is a relative maximal of $E$.
\end{lema}
\pf
Follows the same steps as the proof of \cite[Lemma 1.3]{D87}
\cqd

\begin{lema}\label{absolutecrit} Given a value set $E\subset\mathbb{Z}^d$ and $\alpha\in E$, assume that there exists an index $i\in I$ such that $F_J(E,\alpha)=\emptyset$ for every $J\subsetneq I$ with $i\in J$. Then $\alpha$ is an absolute maximal of $E$.
\end{lema}
\pf
We have to prove that $F_K(E,\alpha)=\emptyset$ for all $K\subset I$ with $i\notin K$. 

Assume, by reductio ad absurdum, that there exists some $K\subset I$ with $i\notin K$ such that  $F_K(E,\alpha)\neq\emptyset$. Let $\beta$ be an element in $F_K(E,\alpha)$, then $\beta_k=\alpha_k$, $\forall k\in K$ and $\beta_j>\alpha_j$, for all $j\notin K$. Applying Property (B) for $\alpha$, $\beta$ and any index $k'\in K$, we have that there exists $\theta\in E$ such that $\theta_{k'}>\beta_{k'}=\alpha_{k'}$, $\theta_l\geq \min\{\alpha_l,\beta_l\}$, $\forall l\neq k'$ and $\theta_j=\alpha_j$ for all $j\notin K$ . If $B=(I\setminus K)\cup\{l\in K, \theta_l=\alpha_l\}$, then we have $\theta\in F_B(E,\alpha) \, (\neq\emptyset)$, with $i\in B$, which is a contradiction.
\cqd

\section{Colengths of fractional ideals}

Let $R$ be a complete admissible ring and let $\J\subset \I$ two regular fractional ideals of $R$ with value sets $D$ and $E$, respectively. Since $\J\subset \I$, one has that $D\subset E$, hence $c(E)\leq c(D)$. Our aim in this section is to find a formula for the length $\ell_R(\I/\J)$ of $\I/\J$ as $R$-modules, called the colength of $\J$ with respect to $\I$, in terms of the value sets $D$ and $E$. 

The motivation comes from the case $r=1$, that is, when $R$ is a domain. In this case, as observed by Gorenstein \cite{Go}, one can easily show that
\[
\ell_R(\I/\J)=\#(E\setminus D).
\]\smallskip
When $r>1$, then $E\setminus D$ is not finite anymore.\medskip

For $\alpha\in \Z^r$ and $\I$ a fractional ideal of $R$, with value set $E$, we define 
\[
\I(\alpha)=\{h\in \I\setminus Z(K); \ v(h)\geq \alpha\}.
\]

It is clear that if $m_E=\min E$, then $\I(m_E)=\I$.\medskip

One has the following result:

\begin{prop}{\rm (\cite[Proposition 2.7]{BDF1})}\label{prop2}
Let $\J\subseteq\I$ be two fractional ideals of $R$, with value sets $D$ and $E$, respectively, then  
$$\ell_R\left(\dfrac{\I}{\J}\right) = \ell_R\left( \dfrac{\I}{\I(\gamma)}\right)-\ell_R\left(\dfrac{\J}{\J(\gamma)}\right),$$ 
for sufficiently large $\gamma\in \N^r$ (for instance, if $\gamma \geq c(D)$).
\end{prop}

If $e_i\in \Z^r$ denotes the vector with zero entries except the $i$-th entry which is equal to $1$, then the following result will give us an effective way to calculate colengths of ideals.

\begin{prop}{\rm \cite[Proposition 2.2]{Da}} \label{prop1}
If $\alpha\in\mathbb{Z}^r$, then we have
$$\ell_R \left(\dfrac{\I(\alpha)}{\I(\alpha+e_i)}\right)=\left\{
\begin{array}{ll}
1, & \ if \  \overline{F}_i(E,\alpha)\neq\emptyset, \\ \\
0, & \ \text{otherwise}.
\end{array}
\right.$$
\end{prop}

So, to compute, for instance, $\ell_R\left(\dfrac{\I}{\I(\gamma)}\right)$, one may take a chain 
$$m_E=\alpha^0 \leq \alpha^1 \leq \cdots \leq \alpha^m=\gamma,$$
where $\alpha^j\in \Z^r$ and $\alpha^{j}-\alpha^{j-1}\in \{e_i, \; i=1,\ldots,r\}$, and then using Proposition \ref{prop1} by observing that 
$$ 
\ell_R\left(\dfrac{\I}{\I(\gamma)}\right)=\ell_R\left(\dfrac{\I(\alpha^0)}{\I(\gamma)}\right)=\sum_{j=1}^m \ell_R\left(\dfrac{\I(\alpha^{j-1})}{\I(\alpha^j)}\right).
$$
D'Anna in \cite{Da} showed that $\ell_R\left( \dfrac{\I}{\I(\gamma)}\right)$ is equal to the length $n$ of a saturated chain $m_E<\alpha^0<\alpha^1<\cdots<\alpha^n=\gamma$ in $E$. The drawback of this result is that one has to know all points of $E$ in the hypercube with opposite vertices $m_E$ and $\gamma$.\medskip

The fact that $E$ is determined by its projections $E_J$ and its relative maximal points, suggests that $\ell_R\left( \dfrac{\I}{\I(\gamma)}\right)$ can be computed in terms of these data. In fact, this will be done in Theorem 1 below.\medskip

In what follows we will denote $\ell_R$  simply by $\ell$.

\subsection{Case r=2}

This simplest case was studied by Barucci,  D'Anna and Fr\"oberg in \cite{BDF1} and we reproduce it here because it gives a clue on how to proceed in general. 

Let $\alpha^0=m_E$ and consider the chain in $\mathbb{Z}^2$ 
$$\alpha^0\leq \cdots\leq\alpha^m=\gamma=(\gamma_1,\gamma_2)\geq c(E)$$
such that
$$
\begin{array}{l}
\alpha^0=(\alpha^0_1,\alpha^0_2),\, \alpha^1=(\alpha^0_1+1,\alpha^0_2),\ldots,\alpha^s=(\gamma_1,\alpha^0_2),\\ \\
\alpha^{s+1}=(\gamma_1,\alpha^0_2+1),\, \alpha^{s+2}=(\gamma_1,\alpha^0_2+2) , \ldots,\alpha^m=(\gamma_1,\gamma_2),
\end{array}$$
and consider the following sets
$$L_1=\{\alpha^0,\alpha^1,\ldots, \alpha^s\} \ \ \text{and} \ \ L_2=\{\alpha^s,\alpha^{s+1},\ldots, \alpha^m\}.
$$

By Proposition \ref{prop1}, we have
\begin{eqnarray}
\ell\left(\dfrac{\I}{\I(\gamma)}\right)&=&\#L_1-\#\{\alpha\in L_1; \ \overline{F}_1(E,\alpha)=\emptyset\}+\nonumber \\
&&\#L_2-\#\{\alpha\in L_2; \ \overline{F}_2(E,\alpha)=\emptyset\}.\nonumber
\end{eqnarray}

Now, because of our choice of $L_1$, denoting by $\mathcal{L}(E_i)$ the set of gaps of $E_i$ in the interval $(\min(E_i), +\infty)$, we have that 
$$
\forall \, \alpha\in L_1, \ \overline{F}_1(E,\alpha)=\emptyset \ \Longleftrightarrow \ pr_1(\alpha)\in \mathcal{L}(E_1),$$
hence 
$$\#\{\alpha\in L_1; \ \overline{F}_1(E,\alpha)=\emptyset\}=\# \mathcal{L}(E_1).$$

Observe that not all $\alpha\in L_2$ with $\overline{F}_2(E,\alpha)=\emptyset$ are such that $pr_2(\alpha)\in \mathcal{L}(E_2)$, hence
$$\#\{\alpha\in L_2; \ \overline{F}_2(E,\alpha)=\emptyset\}=\# \mathcal{L}(E_2)-\xi,$$
where $\xi$ is the number of $\alpha$  in $L_2$  with $pr_2(\alpha)\in E_2$ and $\overline{F}_2(E,\alpha)=\emptyset$. But such $\alpha$ are in one-to-one correspondence with the maximals of $E$, hence $\xi=\#M(E)$.\medskip

Putting all this together, we get
\begin{prop} If $\gamma\geq c(E)$, then 
\begin{equation}\label{I} 
\ell\left(\dfrac{\I}{\I(\gamma)}\right)= (\gamma_1-\alpha^0_1) +(\gamma_2-\alpha^0_2)-\# \mathcal{L}(E_1) -\# \mathcal{L}(E_2) -\#M(E).
																			\end{equation}
\end{prop}

\subsection{Case $r\geq 3$}

Let us assume that $\I$ is a fractional ideal of $R$, where $R$ has $r$ minimal primes.

\begin{figure}[H]
 \centering
\includegraphics[scale=0.3]{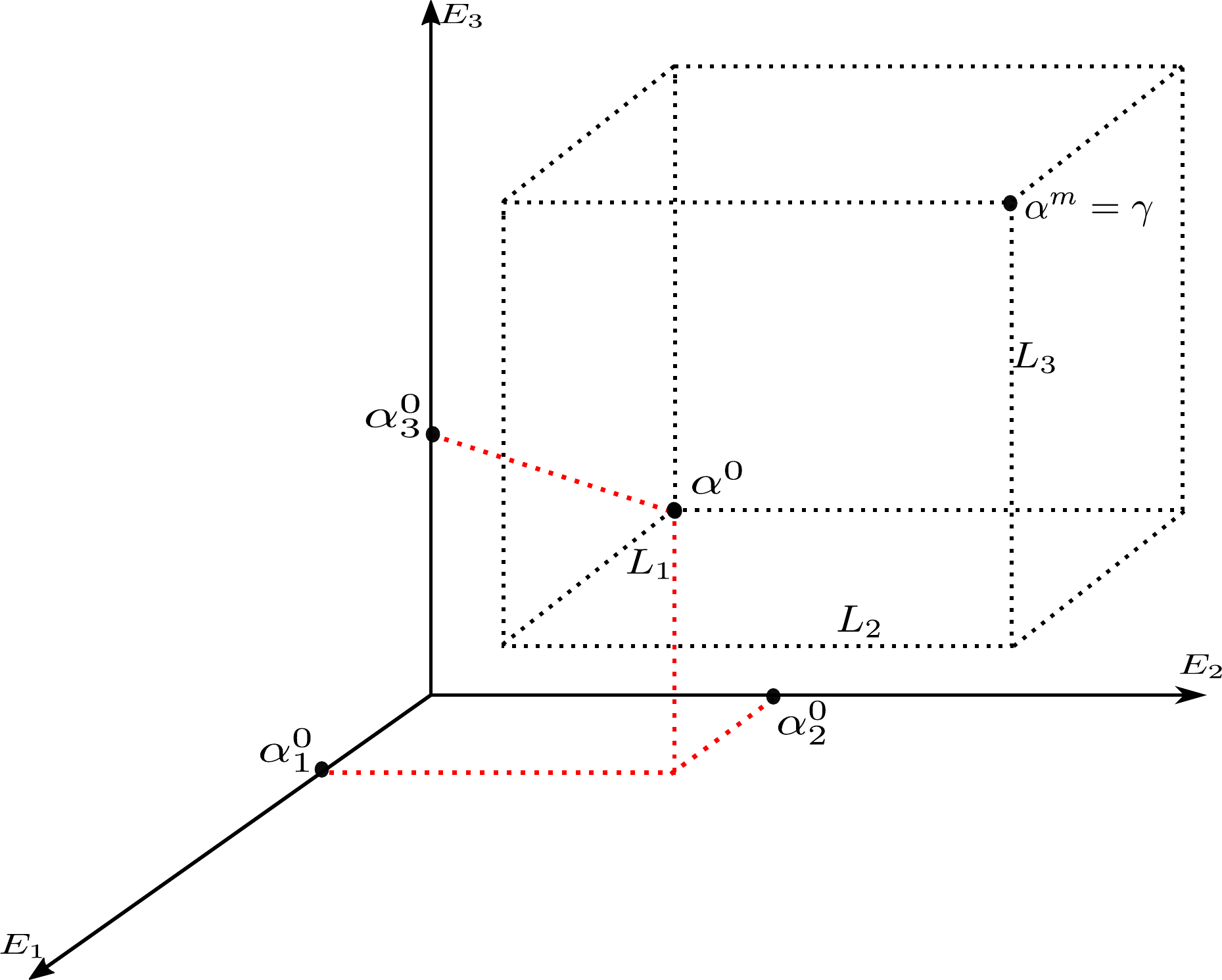}
\caption{The chain for $r=3$}
\label{fig2}
\end{figure}

Let  
 $$m_E=\alpha^0\leq \alpha^1\leq \cdots\leq \alpha^m=\gamma\geq c(E),$$
be the chain in $\mathbb{Z}^r$, given by the union of the following paths (see  Figure \ref{fig2}, for $r=3$):
$$
\begin{array}{l}
L_1\colon \alpha^0, \alpha^1= \alpha^0+e_1,\ldots,\alpha^{s_1}=\alpha^0+(\gamma_1-\alpha^0_1)e_1=(\gamma_1,\alpha^0_2,\ldots,\alpha^0_r),\\ \\
\ldots \\ \\
L_r\colon \alpha^{s_{r-1}}=(\gamma_1,\ldots,\gamma_{r-1},\alpha^0_r), \alpha^{s_{r-1}+1}= \alpha^{s_{r-1}}+e_r, \ldots, \alpha^{m}=
\gamma.
\end{array}
$$

Let us define $I'=\{1,\ldots,r-1\}$. We will need the following result:

\begin{lema}\label{claim} For any $\alpha\in L_1\cup\ldots \cup L_{r-1}$, and for $i\in I'=\{1,\ldots,r-1\}$, one has
$$\overline{F}_i(E,\alpha)\neq\emptyset \ \Longleftrightarrow \ \overline{F}_i(E_{I'},pr_{I'}(\alpha))\neq\emptyset.$$ 
\end{lema}
\pf ($\Rightarrow$) \ This is obvious.

\noindent ($\Leftarrow$)  \ Suppose that $$(\theta_1,\ldots,\theta_{r-1})\in \overline{F}_i(E_{I'}, pr_{I'}(\alpha))\neq\emptyset.$$  Since by Lemma \ref{proj} one has that $pr_{I'}(E)=E_{I'}$, then there exists $\theta=(\theta_1,\ldots,\theta_{r-1},\theta_r)\in E$. Since $\alpha\in L_i$ for some $i=1,\ldots,r-1$, it follows that $\alpha_r=\alpha^0_r$. Then one cannot have $\theta_r<\alpha_r=\alpha^0_r$, because otherwise $$(\alpha^0_1,\ldots,\alpha^0_{r-1},\theta_r)=\min(\alpha^0,\theta)\in E,$$
which is contradiction, since  $\alpha^0$ is the minimum of $E$. Hence $\theta_r\geq\alpha_r$, so $\theta \in \overline{F}_i(E,\alpha)$, and the result follows. 
\cqd

Lemma \ref{claim} allows us to  write:
\begin{equation}\label{intermediateformula}
\ell\left(\dfrac{\I}{\I(\gamma)}\right)=\ell\left(\dfrac{\pi_{I'}(\I)}{\pi_{I'}(\I)(pr_{I'}(\gamma))}\right)+(\gamma_r-\alpha^0_r)-
                                       \#\{\alpha\in L_r; \ \overline{F}_r(E,\alpha)=\emptyset\}.
\end{equation}

Hence to get an inductive formula for $\ell\left(\frac{\I}{\I(\gamma)}\right)$, we only have to compute $$\#\{\alpha\in L_r; \ \overline{F}_r(E,\alpha)=\emptyset\},$$ and for this we will need the following lemma.

\begin{lema}\label{fibra}
Let $\alpha\in\mathbb{Z}^r$, then $\overline{F}_j(E,\alpha)=\emptyset$ if and only if either $\alpha_j\in \mathcal L(E_j)$ or 
there exist some $J\subseteq I$ with $\{j\}\subsetneq J$ and a relative maximal $\beta$ of $E_J$ such that $\widetilde{pr}_{j}(\beta)=\alpha_j$ and $\widetilde{pr}_i(\beta)<\alpha_i$, for all $i\in J$, $i\neq j$.
\end{lema}
\pf ($\Leftarrow$) \ (We prove more, since it is enough to assume $\beta$ is any maximal of $E_J$) It is obvious that if $\alpha_j\in \mathcal L(E_j)$, then $\overline{F}_j(E,\alpha)=\emptyset$. Let us now assume that there exist $ J\subset I$, with $\{j\} \subsetneq J$ and $\beta\in M(E_J)$, such that $\widetilde{pr}_{j}(\beta)=\alpha_j$ and $\widetilde{pr}_i(\beta)<\alpha_i$, for all $i\in J$, $i\neq j$. 

Suppose by reductio ad absurdum that  $\overline{F}_j(E,\alpha)\neq\emptyset$. Let $\theta\in\overline{F}_j(E,\alpha)$ that is $\theta_j=\alpha_j$ and $\theta_i\geq\alpha_i, \forall i\in J\setminus \{j\}$. Now since, $\forall i\in J, i\neq j$, 
$$\widetilde{pr}_j(pr_J(\theta))=\theta_j=\alpha_j=  \widetilde{pr}_{j}(\beta) \ \text{and} \ \widetilde{pr}_i(pr_J(\theta))=\theta_i\geq \alpha_i>\widetilde{pr}_{i}(\beta),$$
then $pr_J(\theta)\in F_j(E_J,\beta)$, which contradicts the assumption that  $\beta\in M(E_J)$.  
	
\noindent ($\Rightarrow$)  \ Since $\overline{F}_j(E,\alpha)=\emptyset$ implies $F_j(E,\alpha)=\emptyset$, the proof follows the same lines as the proof of \cite[Theorem 1.5]{D88}.\cqd

Going back to our main calculation, by Lemma \ref{fibra}, if $\alpha\in L_r$ is such that $\overline{F}_r(E,\alpha)=\emptyset$, then either $\alpha_r\in \mathcal L(E_r)$, or there exist a subset $J$ of  $I=\{1,\ldots,r\}$, with $\{r\}\subsetneq J$, and $\beta\in RM(E_J)$, with $\widetilde{pr}_r(\beta)=\alpha_r$ and $\widetilde{pr}_i(\beta)<\alpha_i$ for $i\in J, i\neq r$. 

Notice that for $\alpha \in L_r$ one has $\alpha_i=\gamma_i$ for $i\neq r$, so the condition $\widetilde{pr}_i(\beta)<\alpha_i$ for $i\in J, \, i\neq r$ is satisfied, since $\beta\in M(E_J)$. So, we have a bijection 
$$ \{\alpha\in L_r; \ \overline{F}_r(E,\alpha)=\emptyset\} \quad \longleftrightarrow  \quad \mathcal L(E_r) \cup \bigcup_{\{r\}\subsetneq  J\subseteq I}\widetilde{pr}_r(RM(E_J)).$$

Since for all $J$,  with $\{r\}\subsetneq J\subseteq I$, the sets $\mathcal{L}(E_r)$ and $\widetilde{pr}_r(RM(E_{J}))$ are disjoint,  it follows that 
\begin{equation} \label{formula4}\#\{\alpha\in L_r;\overline{F}_r(E,\alpha)=\emptyset\}=\# \mathcal{L}(E_r)+\#\bigg(\bigcup_{\{ r\} \subsetneq J\subset I}\widetilde{pr}_r(RM(E_J))\bigg). 
\end{equation}

Now, putting together Equations (\ref{intermediateformula}) and (\ref{formula4}), we get the following recursive formula: 

\begin{teo} \label{teogeral1} Let $\I$ be a fractional ideal of a ring $R$ with $r$ minimal primes with values set $E$. If $\gamma\geq c(E)$, then
\begin{equation} \label{formulafinal1}\begin{array}{rcl}
\ell\left(\dfrac{\I }{\I(\gamma)}\right)=\ell\left(\dfrac{\pi_{I'}(\I)}{\pi_{I'}(\I)(pr_{I'}(\gamma))}\right)&+&(\gamma_r-\alpha^0_r)-\#\mathcal{L}(E_r)-\\
& & \#\bigg(\bigcup_{\{ r\} \subsetneq J\subseteq I}\widetilde{pr}_r(RM(E_J))\bigg).
\end{array}
\end{equation}
\end{teo}

\section{A closed formula for $r=3$}

In this section, we provide a nicer formula than Equation (\ref{formulafinal1}), when $r=3$. To simplify notation, for any $J\subset I=\{1,2,3\}$, we will denote by $RM_J, AM_J$ and $M_J$ the sets $RM(E_J), AM(E_J)$ and $M(E_J)$, respectively. Notice also that if $\#J=2$, then $RM_J=AM_J=M_J$.
\medskip

From Formulas (\ref{I}) and (\ref{intermediateformula}),  for $\gamma\geq c(E)$, one has
\[
\begin{array}{rcl}
\displaystyle \ell\left(\frac{\I}{\I(\gamma)}\right) &=& (\gamma_1-\alpha^0_1)-\# \mathcal{L}(E_1)+(\gamma_2-\alpha^0_2)-\# \mathcal{L}(E_2)-\# M_{\{1,2\}}+\\ \\
&& (\gamma_3-\alpha^0_3)-\# \{\alpha\in L_3; \ \overline{F}_3(E,\alpha)=\emptyset\}.
\end{array}
\]

We will use the following notation:
$$L_3'=\{\alpha\in L_3; \overline{F}_3(E,\alpha)=\emptyset\}.$$

Now, from  Lemma \ref{fibra}, the points $\alpha=(\alpha_1,\alpha_2,\alpha_3)\in L_3'$ are such that $\alpha_3\in \mathcal L(E_3)$ or they are associated to maximal points of either $E_{\{1,3\}}$, $E_{\{2,3\}}$, or $E$ with last coordinate equal to $\alpha_3$. So, we have 
\begin{equation}\label{eta}
\# L_3'=\# \mathcal{L}(E_3)+\#M_{\{1,3\}}+\#M_{\{2,3\}}+\#RM-\eta,
\end{equation}
where $\eta$ is some correcting term which will take into account the eventual multiple counting of maximals having the same last coordinate.

To compute $\eta$ we will analyze in greater detail the geometry of maximal points.

If $\alpha, \beta\in M$ with $\alpha_3=\beta_3$, then $\alpha_1\neq \beta_1$ and $\alpha_2\neq \beta_2$. If $\alpha_1<\beta_1$, then necessarily $\beta_2<\alpha_2$.

We say that two relative (respectively, absolute) maximals $\alpha$ and $\beta $ of  $E$ with $\alpha_3=\beta_3$ and $\alpha_1<\beta_1$ are {\em adjacent}, if there is no $(\theta_1,\theta_2,\alpha_3)$ in $RM$ (respectively, in $AM$) with $\alpha_1<\theta_1<\beta_1$ and $\beta_2<\theta_2<\alpha_2$.\medskip

We will describe below the geometry of the maximal points of $E$

\begin{lema}\label{MA}
If $\alpha\in AM$, then one of the following three conditions is verified:
\begin{enumerate}[\rm (i)]
\item there exist two adjacent relative maximals $\beta$ and $\theta$ of $E$ such that $pr_{\{1,3\}}(\beta)=pr_{\{1,3\}}(\alpha)$ and $pr_{\{2,3\}}(\theta)=pr_{\{2,3\}}(\alpha)$;

\item there exists $\beta\in RM$ such that $pr_{\{1,3\}}(\beta)=pr_{\{1,3\}}(\alpha)$ and $pr_{\{2,3\}}(\alpha)\in M_{\{2,3\}}$, or $pr_{\{2,3\}}(\beta)=pr_{\{2,3\}}(\alpha)$ and $pr_{\{1,3\}}(\alpha)\in M_{\{1,3\}}$;

\item   $pr_{\{1,3\}}(\alpha)\in M_{\{1,3\}}$ and $pr_{\{2,3\}}(\alpha)\in M_{\{2,3\}}$.
\end{enumerate}
\end{lema}
\pf 
Let $\alpha=(\alpha_1, \alpha_2, \alpha_3)\in AM$, then $F(E,\alpha)=\emptyset$. We consider the following sets:

$$R_1=\{\beta\in\mathbb{Z}^3; \ \beta_3=\alpha_3, \beta_1>\alpha_1, \beta_2<\alpha_2\}$$ and $$R_2=\{\theta\in\mathbb{Z}^3; \ \theta_3=\alpha_3, \theta_1<\alpha_1, \theta_2>\alpha_2\}.$$ Then there are four possibilities:
	\begin{align*}
	 R_1\cap E\neq \emptyset \ \text{and} \ R_2\cap E\neq \emptyset, &\quad
	 R_1\cap E\neq \emptyset \ \text{and} \ R_2\cap E = \emptyset.\\ \\
	 R_1\cap E = \emptyset  \ \text{and} \  R_2\cap E \neq \emptyset,&\quad
	 R_1\cap E = \emptyset  \  \text{and}  \  R_2\cap E = \emptyset.
	\end{align*}
	
Suppose $R_1\cap E\neq \emptyset$ and $R_2\cap E\neq \emptyset$. Choose  $\beta\in R_1\cap E$ and $\theta\in R_2\cap E$, such that $\alpha_2-\beta_2$ and $\alpha_1-\theta_1$ are as small as possible. Then by Property (A), we have $\min(\alpha,\beta), \min(\alpha,\theta) \in E$. Obviously $pr_{\{1,3\}}(\beta)=pr_{\{1,3\}}(\alpha)$ and $pr_{\{2,3\}}(\theta)=pr_{\{2,3\}}(\alpha)$.
Moreover, according to Lemma \ref{carrel}, these are relative maximals because $F_3(E,\min(\alpha,\beta))$ and $F_3(E,\min(\alpha,\theta))$ are empty and the sets 
$F_{\{1,3\}}(E,\min(\alpha,\beta))$, $F_{\{1,3\}}(E,\min(\alpha,\theta))$, $F_{\{2,3\}}(E,\min(\alpha,\beta))$ and  $F_{\{2,3\}}(E,\min(\alpha,\theta))$ are nonempty. It follows that $\min(\alpha,\beta)$ and $\min(\alpha,\theta)$ are adjacent relative maximals.
	
Suppose $R_1\cap E\neq \emptyset$ and $R_2\cap E = \emptyset$. Choose $\beta\in R_1\cap E$ such that $\alpha_2-\beta_2$ is as small as possible, then, as we argued above, we have that $\min(\alpha,\beta)\in RM$ and $pr_{\{1,3\}}(\beta)=pr_{\{1,3\}}(\alpha)$. Moreover, as $R_2\cap E = \emptyset$, it follows that $pr_{\{2,3\}}(\alpha)\in M_{\{2,3\}}$.
	
The case $R_1\cap E = \emptyset$ and $R_2\cap E \neq \emptyset$ is similar to the above one, giving us the second possibility in (ii).
	
Suppose $R_1\cap E = \emptyset$ and $R_2\cap E = \emptyset$. It is obvious that $pr_{\{1,3\}}(\alpha)\in M_{\{1,3\}}$ and $pr_{\{2,3\}}(\alpha)\in M_{\{2,3\}}$.
\cqd

Given two points $\theta^1,\theta^2\in \Z^3$ such that ${pr}_3(\theta^1)={pr}_3(\theta^2)$, we will denote by $\R(\theta^1,\theta^2)$ the  parallelogram determined by the coplanar points $\theta^1,\theta^2,\min(\theta^1,\theta^2)$ and $\max(\theta^1,\theta^2)$. We have the following result:
\begin{cor}\label{2maximals} Let $\theta^1,\theta^2\in AM$ be such that ${pr}_3(\theta^1)={pr}_3(\theta^2)$. Then one has $\R(\theta^1,\theta^2)\cap RM\neq \emptyset$.
\end{cor}
\pf Because $\theta^1,\theta^2\in AM$, it follows immediately that (iii) of Lemma \ref{MA} cannot happen, therefore, the existence of the relative maximal is ensured by (i) or (ii). \cqd

\begin{lema}\label{adjacent}
If $\beta$ and $\beta'$ are adjacent relative maximals, with $\beta_3=\beta_3'$, then $\max(\beta, \beta')$ is an absolute maximal of $E$.
\end{lema}
\pf We may suppose that $\beta_1>\beta_1'$ and $\beta_2<\beta_2'$. As $\beta$ and $\beta'$ are adjacent, we have that $F_{\{1,3\}}(E,\beta)\cap F_{\{2,3\}}(E,\beta')\neq\emptyset$, because otherwise, take $\alpha^1\in F_{\{1,3\}}(E,\beta)$, with $\alpha^1_2$ the greatest possible and $\alpha^2\in F_{\{2,3\}}(E,\beta')$, with $\alpha^2_1$ the greatest possible. From Lemma \ref{absolutecrit} it follows that $\alpha^1$ and $\alpha^2$ are absolute maximals of $E$, then by Corollary \ref{2maximals} there exists a relative maximal in the region ${\R}(\alpha^1,\alpha^2)$, this contradicts the fact that $\beta$ and $\beta'$ are adjacent relative maximals.

Then, effectively, $F_{\{1,3\}}(E,\beta)\cap F_{\{2,3\}}(E,\beta')=\{\max(\beta,\beta')\}$, which is an absolute maximal.
\cqd

Recall that the elements in $L_3'$ are of the form $(\gamma_1,\gamma_2,\alpha_3)$, with $\alpha^0_3\leq \alpha_3 \leq \gamma_3$.

\begin{lema}\label{M13MR}
Let $\alpha\in L_3'$ be such that $\alpha_3\in (\widetilde{pr}_3(M_{\{1,3\}})\setminus \widetilde{pr}_3(M_{\{2,3\}})) \cap {pr}_3(RM)$, or $\alpha_3\in (\widetilde{pr}_3(M_{\{2,3\}})\setminus \widetilde{pr}_3(M_{\{1,3\}}))\cap {pr}_3(RM)$, then there are the same number of relative as absolute maximals in $E$ with third coordinate equal to $\alpha_3$.
\end{lema}
\pf 
We assume that $\alpha_3\in (\widetilde{pr}_3(M_{\{1,3\}})\setminus \widetilde{pr}_3(M_{\{2,3\}}))\cap {pr}_3(RM)$, since the other case is analogous.

Since $\alpha_3\in {pr}_3(RM)$, we may assume that there are $s \,(\geq 1)$ relative maximals $\beta^1,\ldots,\beta^s$ in $E$ with third coordinate equal to $\alpha_3$. We may suppose that $\beta_1^1<\beta_1^2<\cdots < \beta_1^s$, so the $\beta^i$'s are successively adjacent relative maximals, hence, by lemma \ref{adjacent}, we have that $$\max(\beta^1,\beta^2), \ldots,\max(\beta^{s-1},\beta^s)\in AM.$$ 

This shows that there are at least $s-1$ absolute maximals in $E$ with third coordinate $\alpha_3$.  

Now as ${pr}_3(\alpha)\in \widetilde{pr}_3(M_{\{1,3\}})$, then there is a $(\eta^1_1,\alpha_3)\in M_{\{1,3\}}$ with $\eta^1_1\leq\alpha_1 \,(=\gamma_1)$, because $c(E_{\{1,3\}})\leq pr_{\{1,3\}}(c(E))=(\gamma_1,\gamma_3)$. Because of our hypothesis, the elements $\delta$  in the fiber $F_{\{1,3\}}(E,\beta^s)$ are such that  $\beta^s_1<\delta_1\leq\eta^1_1$. But we must have $\delta_1= \eta^1_1$, because, otherwise, there would be a point $\eta^1=(\eta^1_1,\eta^1_2,\alpha_3)\in pr_{\{1,3\}}^{-1}(\eta^1_1,\alpha_3)$, with $\eta^1_2<\beta^s_2$, and a point $\eta^2\in F_{\{2,3\}}(E,\beta^s)$ with $\eta^2_1<\eta^1_1$ and $\eta^2_2=\beta^s_2$. These $\eta^1$ and $\eta^2$ are absolute maximals, due to Lemma \ref{absolutecrit}, then from Corollary \ref{2maximals}, there would exist a relative maximal in the region $\R(\eta^1,\eta^2)$, which contradicts the fact that we have $s$ relative maximals. This implies that $(\beta^s_1,\eta^1_2,\alpha_3)$ is an absolute maximal of $E$. 

We have to show that there are no other absolute maximals. If such maximal existed, then one of the three conditions in Lemma \ref{MA} would be satisfied. Obviously conditions (i) and (iii) cannot be satisfied, but neither condition (ii) can be satisfied, because otherwise $\alpha_3\in \widetilde{pr}_3(M_{\{2,3\}})$, which is a  contradiction.
\cqd

\begin{lema}\label{M13M23}
Let $\alpha=(\alpha_1,\alpha_2,\alpha_3)\in L_3'$ be such that $\alpha_3\in \big(\widetilde{pr}_3(M_{\{1,3\}})\cap \widetilde{pr}_3(M_{\{2,3\}})\big)\setminus {pr}_3(RM)$, then there exists one and only one absolute maximal of $E$ with third coordinate equal to $\alpha_3$.
\end{lema}	
\pf 
As $\alpha_3\in \widetilde{pr}_3(M_{\{1,3\}})\cap \widetilde{pr}_3(M_{\{2,3\}})$, then there exist $(\beta^1_1,\alpha_3)\in M_{\{1,3\}}$ and $(\beta^2_2,\alpha_3)\in M_{\{2,3\}}$ such that $\beta^1_1<\alpha_1(=\gamma_1)$ and $\beta^2_2<\alpha_2(=\gamma_2)$, because one always has that $c(E_{\{i,j\}}) \leq pr_{\{i,j\}}(c(E))$.

Consider the element $\theta=(\beta_1^1,\beta_2^2,\alpha_3)$. If $\theta\in E$, since it is easy to verify that $F_J(E,\theta)=\emptyset$ for $3\in J\subsetneq \{1,2,3\}$, it follows by Lemma \ref{absolutecrit} that $\theta$ is an absolute maximal of $E$, which is unique in view of Corollary \ref{2maximals} and the hypothesis that $\alpha_3\not\in {pr}_3(RM)$.

If $\theta \not\in E$, then take $\theta_1=(\beta^1_1,\delta^1_2,\alpha_3)\in pr^{-1}_{\{1,3\}}(\beta^1_1,\alpha_3)\cap E$, and $\theta_2=(\delta^2_1,\beta^2_2,\alpha_3)\in pr^{-1}_{\{2,3\}}(\beta^2_2,\alpha_3)\cap E$. We have that $\delta^2_1 <\beta^1_1$ and $\delta^1_2< \beta^2_2$, because otherwise $\theta\in E$ or, $(\beta^1_1,\alpha_3)$ and/or $(\beta^2_2,\alpha_3)$ would not be maximals of $E_{\{1,3\}}$ and/or $E_{\{2,3\}}$. Choose $\delta^2_1$ and $\delta^1_2$ the greatest possible, then it is easy to verify that $F_J(E,\theta_i)=\emptyset$ for $i=1,2$ and $3\in J\subsetneq \{1,2,3\}$. Hence from Lemma \ref{absolutecrit}, $\theta_1$ and $\theta_2$ are absolute maximals of $E$, therefore from Corollary \ref{2maximals} there would be a relative maximal of $E$ with third coordinate equal to $\alpha_3$, which is  a contradiction.
\cqd

\begin{lema}\label{MMMR}
Let $\alpha\in L_3'$ be such that $\alpha_3\in \widetilde{pr}_3(M_{\{1,3\}})\cap \widetilde{pr}_3(M_{\{2,3\}})\cap {pr}_3(RM)$. If there exist $s$ relative maximals with third coodinate equal to $\alpha_3$, then there exist $s+1$ absolute maximals with third coordinate equal to $\alpha_3$.
\end{lema}	
\pf 
Following the proof of Lemma \ref{M13MR}, we have $s-1$ absolute maximals obtained by taking the maximum of each pair of adjacent relative maximals. The conditions $\alpha_3\in \widetilde{pr}_3(M_{\{1,3\}})$ and $\alpha_3\in \widetilde{pr}_3(M_{\{2,3\}})$ give us two extra absolute maximals, and the same argument used there, shows that there are no other.
\cqd

\begin{lema}\label{MR}
Let $\alpha\in L'_3$ be such that $\alpha_3\in {pr}_3(RM)\setminus \big(\widetilde{pr}_3(M_{\{1,3\}})\cup \widetilde{pr}_3(M_{\{2,3\}})\big)$. If there exist $s$ relative maximals with third coordinate equal to $\alpha_3$, then we have $s-1$ absolute maximals with third coordinate equal to $\alpha_3$.
\end{lema}	
\pf The arguments used in the proofs of the last two lemmas give us the result.
\cqd

Going back to Formula (\ref{eta}), we want to calculate $\eta$.  From Lemma \ref{fibra} we can ensure that $\alpha\in L'_3=\{\alpha\in L_3; \ \overline{F}_3(E,\alpha)=\emptyset\}\setminus {\mathcal L}(E_3)$, only if $\alpha$ falls into one of the following five cases:

\begin{enumerate}[\rm (i)]
\item  $\alpha_3\in (\widetilde{pr}_3(M_{\{1,3\}})\setminus \widetilde{pr}_3(M_{\{2,3\}}))\cap {pr}_3(RM)$. 

If there exist such $\alpha$, then they are related to a unique element of $M_{\{1,3\}}$ and if there are $s_1$ relative maximals with third coordinate $\alpha_3$, then in our formula $\alpha$ was counted $s_1+1$ times. By Lemma \ref{M13MR} we know that there exist $s_1$ absolute maximals of $E$ with third coordinate $\alpha_3$. So, we subtract $s_1$ from our counting to partially correct the formula.  

\item  $\alpha_3\in (\widetilde{pr}_3(M_{\{2,3\}})\setminus \widetilde{pr}_3(M_{\{1,3\}}))\cap {pr}_3(RM)$.

Analogously to (i), $\alpha$ is related to a unique element of $M_{\{2,3\}}$ and if there are $s_2$ relative maximals with third coordinate $\alpha_3$, then $\alpha$ was counted $s_2+1$ times in the formula. Again, by Lemma \ref{M13MR} we know that there are $s_2$ absolute maximals of $E$ with third coordinate $\alpha_3$. So, we subtract $s_2$ from our counting to partially correct the formula.  

\item $\alpha_3\in \big(\widetilde{pr}_3(M_{\{1,3\}})\cap \widetilde{pr}_3(M_{\{2,3\}})\big)\setminus {pr}_3(RM)$.

In this case, $\alpha$ is related to a unique elements in $M_{\{1,3\}}$ and in $M_{\{2,3\}}$, so in the formula we are counting  $\alpha$ twice. By Lemma \ref{M13M23} there is a unique absolute maximal of $E$ with third coordinate $\alpha_3$ such that its projections $pr_{\{1,3\}}$ and $pr_{\{2,3\}}$ are in $M_{\{1,3\}}$ and $M_{\{1,3\}}$, respectively. So, we correct partially the formula by subtracting $1$, which corresponds to this unique absolute maximal.

\item $\alpha_3\in \widetilde{pr}_3(M_{\{1,3\}})\cap \widetilde{pr}_3(M_{\{2,3\}})\cap {pr}_3(RM)$.

In this case, $\alpha$ is related to a unique element of $M_{\{1,3\}}$, to a unique element of $M_{\{2,3\}}$ and, let us say, $s_3$ elements of $RM$, so in our counting, $\alpha$ was counted  $s_3+2$ times. By Lemma \ref{MMMR} there exist $s_3+1$ absolute maximals of $E$ with third coordinate $\alpha_3$. In this case, the correcting term is $s_3+1$, equal to the number of these absolute maximals.

\item  $\alpha_3\in {pr}_3(RM)\setminus \big(\widetilde{pr}_3(M_{\{1,3\}})\cup \widetilde{pr}_3(M_{\{2,3\}})\big)$. 

In this case, $\alpha$ is related with, let us say, $s_4$ elements of $RM$ with third coordinate equal to $\alpha_3$, so we are counting it $s_4$ times. By Lemma \ref{MR} there exist $s_4-1$ absolute maximals with third coordinate $\alpha_3$. This is exactly the correcting term we must apply to our formula. 

\end{enumerate}

Observe that the above cases exhaust all absolute maximals of $E$, implying the following result conjectured by M. E. Hernandes after having analyzed several examples (cf. \cite{He}):

\begin{teo}\label{formular=3} Let $R$ be an admissible ring with three minimal primes and let $\I$ be a fractional ideal of $R$ with values set $E$. If $\gamma \geq c(E)$, then 
\[
\begin{array}{rcl}
\ell\left(\frac{\I}{\I(\gamma)}\right)&= &\sum_{i=1}^r \left((\gamma_i-\alpha^0_i)  - \# \mathcal{L}(E_i)\right)  -
																		\sum_{1\leq i<j\leq 3}\#M_{\{i,j\}}-\\ 
																		&&\#RM+\#AM.
\end{array}
\]
\end{teo}

\begin{cor} Let $\J\subseteq\I$ be two fractional ideals of an admissible ring $R$, with three minimal primes. Denote by $E$ and $D$, respectively, the value sets of $\I$ and $\J$. Then  
\[\begin{array}{rcl}
\ell_R\left(\dfrac{\I}{\J}\right) &= &\sum_{i=1}^3\left( (\beta^0_i-\alpha^0_i) + (\#\mathcal L(D_i)- \#\mathcal L(E_i))\right)+ \\                                     && \sum_{1\leq i<j\leq 3}\#M_{\{i,j\}}(D)- \sum_{1\leq i<j\leq 3}\#M_{\{i,j\}}(E)  + \\ \\
																	  &&\# RM(D)-\# RM(E) + \# AM(E)-\# AM(D),
																	\end{array}
																	\]
where $\alpha^0=\min(E)$ and $\beta^0=\min(D)$.
																	\end{cor}

\noindent Author's e-mail addresses: e.marcavillaca@gmail.com \, and \, ahefez@id.uff.br

\bigskip





\begin{thebibliography}{X}

\bibitem{BDF1} \textsc{Barucci, V.}; \textsc{D'anna, M.};  \textsc{Fr\"oberg, R.}, \textit{The Semigroup of Values of a One-dimensional Local Ring with two Minimal Primes}, Communications in Algebra, 28(8), pp 3607-3633 (2000).

\bibitem{Da} \textsc{D'anna, M.}, \textit{The Canonical Module of a One-dimensional Reduced Local Ring}, Communications in Algebra, 25, pp 2939-2965 (1997).

\bibitem{D87} \textsc{Delgado de la Mata, F.}, \textit{The semigroup of values of a curve singularity with several branches}, Manuscripta Math. 59, pp 347-374 (1987).

\bibitem{D88} \textsc{Delgado de la Mata, F.}, \textit{Gorenstein curves and symmetry of the semigroup of values}, Manuscripta Math. 61, pp 285-296 (1988).

\bibitem{Ga} \textsc{Garcia, A.}, \textit{Semigroups associated to singular points of plane curves},
J. Reine. Angew. Math. 336, 165-184 (1982).

\bibitem{Go} \textsc{Gorenstein, D.}, \textit{An arithmetic theory of adjoint plane curves}, Trans. Amer. Math. Soc. 72, pp 414-436 (1952).

\bibitem{He} \textsc{Hernandes, M. E.},  \textit{Private communication}

\bibitem{KST} \textsc{Korell, P.; Schulze, M.; Tozzo, L.}, \textit{Duality on value semigroups}, arXiv:1510.04072v4 [math.AG] (2 Nov. 2017); to appear in J.  Commut. Alg.


\end{thebibliography}
\end{document}